\documentclass[12pt]{amsart}
\usepackage{amsmath,amsfonts,amsthm,amssymb,pstricks,pstcol,pst-plot}
\begin{document}

\def\N{{\mathbb N}}
\def\R{{\mathbb R}}
\def\Z{{\mathbb Z}}
\newcommand{\ds}{\displaystyle}
 \newtheorem{theorem}{Theorem}
 \newtheorem{remark}{Remark}[section]
 \newtheorem{lemma}[remark]{Lemma}
 \newtheorem{cor}[remark]{Corollary}
 \newtheorem{prop}[remark]{Proposition}
 \newtheorem{definition}[remark]{Definition}
\def\etf{{e^{{\tau} \phi}}}
\def\tnf{\tau \nabla \phi}
\def\f12{{\frac{1}{2}}}
\def\ty{\tilde{y}}
\def\rar{\rightarrow}
\def\e{\varepsilon}
\def\g{\gamma}
\def\rt{\tilde{r}}

\def\bb{\beta}
\def\dd{\delta}
\def\d{\partial}
\def\aa{\alpha}
\def\l{\lambda}
\def\th{\tilde{H}}
\def\SUCP{\mbox{\bf (SUCP) }}
\def\dt{\tilde{\Delta}}
\def\nt{\tilde{\nabla}}
\def\pt{\tilde{P}}
\def\td{\tilde{\d}}
\def\tX{\tilde{X}}
\def\tY{\tilde{Y}}
\def\dt{\tilde{P}_0}
\def\oops{{\bf !!!!!}}

\title{ $L^p$ eigenfunction bounds for the Hermite operator}

\author{Herbert Koch}

\address{Fachbereich Mathematik,Universit\"at Dortmund,  44221
  Dortmund}
\email{koch@mathematik.uni-dortmund.de}

\author{Daniel Tataru}
\address{Department of Mathematics, University of California at
  Berkeley, Berkeley, CA, 94596}
\email{tataru@math.berkeley.edu}
\thanks{The second author was supported in part by NSF Grant DMS-0226105}


\begin{abstract}
We obtain $L^p$ eigenfunction bounds for  the harmonic oscillator $H =
-\Delta + x^2$ in $\R^n$ and  for other related operators,
improving earlier results of Thangavelu and Karadzhov. We also
construct suitable counterexamples which show that our
estimates are sharp.

\end{abstract}

\maketitle

\section{Introduction}
The question of obtaining $L^p$ eigenfunction bounds for elliptic
operators on compact manifolds has been considered in Sogge's work,
for which we refer the reader to his book \cite{MR94c:35178}. The
$L^p$ eigenfunction bounds in \cite{MR94c:35178} are sharp, and turn
out to be related to a variable coefficient version of the restriction
theorem, and further to a phase curvature condition for Fourier
integral operators. In this analysis a special role is played by the
Laplace-Beltrami operator for the sphere, which is the worst case because
it has many highly concentrated eigenfunctions. This is connected to
the fact that it has a periodic Hamilton flow.

In this article we consider the problem of obtaining $L^p$ eigenfunction
bounds for the Hermite operator $H = -\Delta + x^2$ in $\R^n$ 
and  also for a larger class of related operators of the form $H_V =  -\Delta +
V$.  Within this class the Hermite operator plays a role similar
to that of the spherical Laplacian, in that it has a periodic Hamilton
flow and many highly concentrated eigenfunctions.

This question has received considerable interest in the context 
of Riesz summability for the harmonic oscillator in the work of
Thangavelu~\cite{MR91b:42048}, \cite{MR94i:42001},
\cite{MR2000a:42044} and Karadzhov~\cite{MR96b:46043}.

Our interest in it has a different source, namely the strong unique
continuation problem for parabolic equations. In this context the work
of Thangavelu and Karadzhov has already found applications in
Escauriaza~\cite{MR2001m:35135} and
Escauriaza-Vega~\cite{MR2003b:35088}.  Further applications are
contained in a forthcoming paper of the authors.  We note in passing
that the related strong unique continuation problem for second order
elliptic operators is related to the eigenfunction bounds for the
spherical harmonics. This was first observed in work of
Jerison~\cite{MR88b:35218}; see also the authors paper
\cite{MR2001m:35075} and further references therein.

In the next section we begin by considering the problem of obtaining
dispersive and Strichartz estimates for the corresponding
Schr\"odinger equation.  This leads to an alternative proof of the
eigenfunction bounds of Karadzhov~\cite{MR96b:46043} and Thangavelu
\cite{MR2000a:42044}. Our approach has the advantage that is robust
enough so that it allows us to obtain the same results with $x^2$
replaced by potentials in a very large class.

Then we direct our attention to the eigenfunctions $(-\Delta + x^2) \phi
= \lambda^2 \phi$. These are concentrated inside the ball $\{|x| \leq
\lambda\}$, and have an exponential Airy type decay beyond this
threshold.  The behavior of the eigenfunctions inside the ball is not
very different from (a rescaling of) what happens in a bounded domain.
However, considerable care is required near the boundary of the ball,
where the concentration scales are different. Consequently, it is more
natural to obtain weighted $L^p$ estimates with weights which are
essentially powers of $\lambda - |x|$. The results we obtain
strengthen those of Karadzhov and Thangavelu and complete the picture.
As before, our methods are robust and apply equally to any
potential which behaves roughly like $x^2$.

In the last section of the paper we construct appropriate examples
which illustrate the possible concentration scales for eigenvalues of
the Hermite operator and show that our $L^p$ bounds are sharp.

To conclude the introduction we provide the reader with a special case
of our main result.  Denote the spectral projection to the
eigenvalue $k=\lambda^2$ by $P_k$.  Then our $L^p$ bounds for
eigenfunctions of the Hermite operator $H$ have the form
\begin{equation} 
\Vert P_k \phi  \Vert_{L^p} \lesssim k^{\rho(p)/2} \Vert \phi
\Vert_{L^2} 
\end{equation} 
where $\rho$ is given as in the following figure.

\begin{figure}[hb] \hspace*{-.8in}
\begin{pspicture}(-.5,-.5)(10,7)
\psset{xunit=20cm,yunit=3cm}
\psaxes(0,.5)(0,0)(.6,2.2)
\psline(0,2)(.3,.5)
\psline(.3,.5)(.375,.375) 
\psline(.375,.375)(.5,.5) 
\rput(-.03,2){$\frac{n-2}2$}
\psline(-.01,2)(.01,2)
\rput(-.03,2.2){$\rho$}
\rput(-.035,.33){$-\frac1{n+3}$}
\psline(-.01,.375)(.01,.375)
\rput(.6,.25){$\frac1p$}
\rput(.5,.25){$\frac12$}
\psline(.5,.45)(.5,.55)
\rput(.3,.25){$\frac{n-2}{2n}$}
\psline(.3,.45)(.3,.55)
\rput(.333,.75){$\frac{n-1}{2(n+1)}$}
\psline(.333,.45)(.333,.55)
\rput(.375,.25){$\frac{n+3}{2(n+1)}$}
\psline(.375,.45)(.375,.55)
\end{pspicture}
\caption{The exponent $\rho$ as function of $1/p$}
\label{fig}\end{figure}
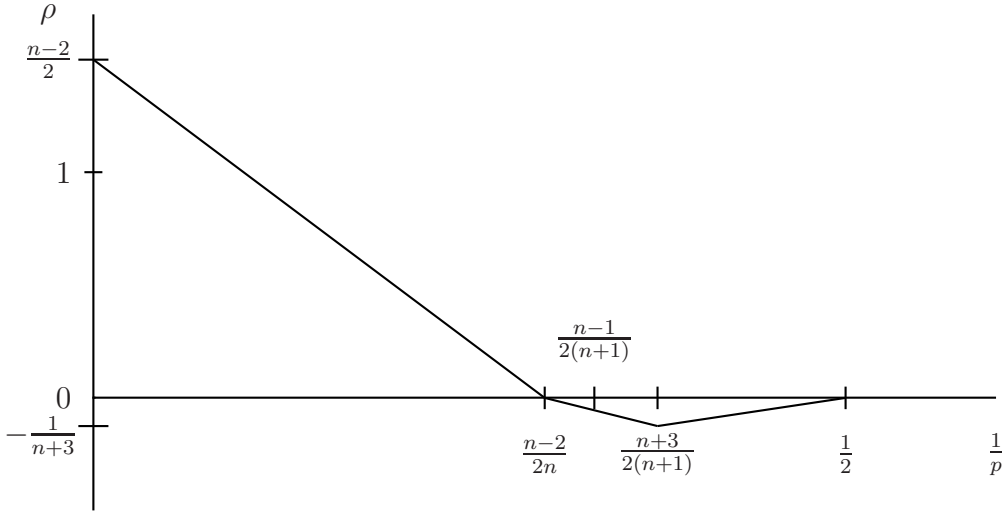

\section{Dispersive estimates for the Schr\"odinger equation}

Consider the  Schr\"odinger operator associated to the
Hermite operator in $\R^n$,
\[
i \d_t - H
\]
which generates a group of isometries $t \to e^{itH}$.
Furthermore, this group is periodic with period $\pi$.
Here we investigate the question of obtaining pointwise
bounds for the kernel of $e^{itH}$. We also consider the
same question for more general operators
\[
H_V = -\Delta + V
\]

\begin{theorem}
a) The operators $e^{itH}$ satisfy
\[
\|e^{itH}\|_{L^1 \to L^\infty} \lesssim |\sin{t}|^{-\frac{n}2}
\]
\label{decay}

b) Let $V$ be a potential which satisfies
\begin{equation}
  \label{vi}
|\d^\alpha V| \leq c_\alpha, \qquad |\alpha| \geq 2  
\end{equation}
Then for  $|t| \ll 1$ the operators $e^{itH}$ satisfy
\[
\|e^{itH_V}\|_{L^1 \to L^\infty} \lesssim |t|^{-\frac{n}2}
\]
\label{decayv}\end{theorem}

\begin{proof}[Proof of Theorem~\ref{decay}]
  Prove the bound at time $t_0$.  In case (a) the estimate follows
  from well known formulas for the Schr\"odinger kernel of the Hermite
  operator. We prefer to give a more flexible argument.  By
  periodicity we can assume that $|t_0| \leq \pi/2$ and replace $|\sin
  t|$ by $|t|$. We rescale to time $1$ by setting
\[
t \to \frac{t}{t_0}, \qquad x \to \frac{x}{\sqrt{t_0}}
\]
Then we need to prove an uniform bound
for the kernel of $e^{iH^{t_0}}$, respectively $e^{iH^{t_0}_V}$ where
\[
H^{t_0}= t_0^2  x^2 -\Delta, \qquad  H^{t_0}_V =  t_0 V(x \sqrt{t_0})-\Delta  
\]
Since $t_0$ is in a bounded set, it follows that the symbols
\[
h^{t_0}(x,\xi) = t_0^2 x^2 +\xi^2, \qquad h^{t_0}_V (x,\xi) = t_0 V(x
\sqrt{t_0})+\xi^2
\]
satisfy the bounds
\begin{equation}
|\d_x^\alpha \d_\xi^\beta h(x,\xi)| \leq c_{\alpha,\beta},
\qquad |\alpha|+|\beta| \geq 2
\label{bo}\end{equation}
uniformly in $|t_0| \lesssim 1$. This implies that we can use directly
Proposition~4.3  in \cite{oi} to obtain a phase space representation of
the fundamental solution $K(t,y,\tilde y)$ for $i \d_t - H^{t_0}_V$. 
For the reader's convenience we restate the result here:

\begin{lemma}[\cite{oi}] Let $h$ be a symbol which satisfies
  \eqref{bo}. Then for $|t| \lesssim 1$ the fundamental solution
$K(t,y,\tilde y)$ for $i \d_t - h^w(x,D)$ can be represented as
\begin{equation}
K(t, y,\tilde y)  = 
  \int_{\mathbb{R}^{2n}}
e^{- \frac12 (\tilde y-x)^2} e^{ 
- i \xi (\tilde y-x)}    e^{i \psi(t,x,\xi)}  e^{i \xi^t( y-x^t)} G(t,x,\xi,y)   dx \, d\xi   
\label{fs} \end{equation}
where $(x,\xi)$ to $(x^t,\xi^t)$ is the Hamilton flow for $h$,
the function $G$ satisfies  
\begin{equation}
|(x^t-y)^\gamma \d_x^\aa \d_\xi^\bb \d_{y}^\nu G(t,x,\xi,y)| \lesssim
c_{\gamma,\aa,\bb,\nu}
\label{Gbd}\end{equation}
and the real phase function $\psi$ is determined by
\[
 \frac{d}{dt} \psi(t,x,\xi) = -h(x^t,\xi^t) + \xi^t h_\xi(x^t,\xi^t),  \qquad \psi(0,x,\xi) = 0
\]
\end{lemma}

A feature of the construction in \cite{oi} is that the integrand
solves the evolution equation for each $(x,\xi)$. In addition, it is
concentrated in the phase space on the unit scale along the
bicharacteristic $t \to (x^t,\xi^t)$. Such highly localized solutions
are called wave packets. With this terminology, one can view the above
lemma as a way of representing solutions for $i \d_t - h^w(x,D)$ as
almost orthogonal superpositions of wave packets.

In our case we need a pointwise bound for $K(1,y,\tilde y)$.
Neglecting all oscillations in \eqref{fs} we write
\[
|K(1,y,\tilde y)| \lesssim  \int_{\mathbb{R}^{2n}} e^{- \frac12 (\tilde y-x)^2}
(1+|y-x^1|)^{-N} dx \, d\xi   
\]
Due to \eqref{bo}  the Hamilton flow  for $H^{t_0}_V$ is Lipschitz. Hence the 
integration in $x$ is trivial, and we obtain
\[
|K(1,y,\tilde y)| \lesssim  \int_{\mathbb{R}^{n}}
(1+|y-{\tilde y}^1|)^{-N}  d\xi   
\]
In order to obtain an uniform bound it suffices to prove that the
Lipschitz map $\xi \to {\tilde y}^1$ has a Lipschitz inverse.  In the
case of the Hermite operator this map is linear, and the desired
conclusion is obtained by direct computation. For a more general
potential $V$ one computes the linearization of the Hamilton flow,
which shows that
\[
\frac{{d \tilde y}^1}{d \xi} = 2 I_n + O(t_0)
\]
For small $t_0$ this shows that the map $\xi \to {\tilde y}^1$
is a global diffeomorphism of $\R^n$.
\end{proof}

As a consequence of the dispersive estimates for the
Schr\"odinger equation we also obtain Strichartz estimates for the
Schr\"odinger equation:

\begin{theorem}
Let $V$ be a potential which satisfies 
\begin{equation}
 \label{v2}
|\d^\alpha V| \lesssim 1, \qquad |\alpha| = 2  
\end{equation}
Then the solution $u$   to 
\[
( i \d_t -H_V ) u = f, \qquad u(0) = u_0
\]
satisfies
\[
\| u \|_{L^{p_1}(0,1;L^{q_1})} \lesssim \|u_0\|_{L^2} + \|f\|_{ L^{p'_2}(0,1;L^{q'_2})}
\]
whenever the pairs $(p_1, q_1)$ and $(p_2,q_2)$ are subject to
\[ 
\frac{2}p + \frac{n-1}{q} = \frac{n-1}2, \qquad 2 \leq p \leq \infty, 
\quad (n,p,q) \neq (2,2,\infty)
\]
\label{se}\end{theorem}

\begin{proof}
  a) If $V$ satisfies the stronger bound \eqref{vi} then this follows
  from Theorem~\ref{decayv} by standard arguments as in
  Ginibre-Velo~\cite{GV} and the references therein. For the endpoint
  one can use the results in Keel-Tao~\cite{KT}.

b)  For potentials $V$ which only satisfy \eqref{v2} we use 
a frequency decomposition of $V$. Given a smooth compactly supported
function $\chi$ which equals $1$ in a neighborhood of the origin we set
\[
V = V_0 + V_1, \qquad V_0 = \chi(D) V, \quad V_1 = (1-\chi(D) V)
\]
 
The low frequency part $V_0$ satisfies \eqref{vi}, therefore we use
part (a) of the proof. On the other hand the high frequency part $V_1$
is bounded, so we can add it in using the fact that the result in the
theorem is stable with respect to $L^2$ bounded perturbations of
$H_V$.
\end{proof}

Using the Strichartz estimates one can easily obtain eigenfunction
bounds. We begin with the Hermite operator $H$. Its spectrum is
$n+2\N$.  A basis of orthogonal eigenfunction is given by the Hermite
functions,
\[ h_\alpha = (-1)^{|\alpha|} e^{x^2/2} \d^\alpha e^{-x^2} \]
with 
\[ H h_\alpha = (n + 2|\alpha|) h_\alpha \]

For an integer $k= \lambda^2$ in the
spectrum we denote by $P_k$ the corresponding projector.
For $\phi \in L^2$ the function 
\[
u = e^{itk} P_k \phi
\]
solves $(i \d_t -H) u =0$. Hence we can apply Theorem~\ref{se} to
obtain:
\begin{cor}
Let $P_k$ be a spectral projector for $H$. Then 
\[
\|P_k \phi\|_{L^p} \lesssim \|\phi\|_{L^2}, \qquad 2 \leq p \leq
\frac{2n}{n-2}, \qquad (n,p) \neq (2,\infty)
\]
\label{lpef}\end{cor}
In particular this gives $L^p$ bounds for each eigenfunction.
The estimate still holds in the case $ (n,p) = (2,\infty)$, but to
prove this one needs the stronger results in the next section.

In the case of the operator $H_V$ the spectrum need not be discrete.
Even if it is, there is no guarantee that there is an $O(1)$ spectral
gap. Hence it is more useful to consider spectral projectors $P_{[k,k+1]}$
associated to unit size intervals. For $\phi \in L^2$ the function 
\[
u = e^{itk} P_{[k,k+1]} \phi
\]
satisfies
\[
\|(i \d_t -H) u\|_{L^\infty L^2} \leq \|\phi\|_{L^2}
\]
 Hence we can apply Theorem~\ref{se} to
obtain:
\begin{cor}
  Let $V$ be a potential which satisfies \eqref{v2}.  If
  $P_{[k,k+1]}$ is a spectral projector for $H_V$ then
\[
\|P_{[k,k+1]} \phi\|_{L^p} \lesssim \|\phi\|_{L^2}, \qquad 2 \leq p \leq
\frac{2n}{n-2}, \qquad (n,p) \neq (2,\infty)
\]
\label{lpefv}\end{cor}

Here and in the sequel we formulate results for $n\ge 2$, which are valid for 
$n=1$ with modifications which are either obvious or discussed. 

\section{Weighted $L^p$ eigenfunction bounds}

Let $(\phi, \l^2)$ be an eigenfunction, respectively an eigenvalue for
the Hermite operator. The function $\phi$ is essentially concentrated
in the ball $\{|x| \leq \l\}$ (modulo an exponentially decaying tail).
We split the interior of this region into overlapping dyadic parts
with respect to the distance to the boundary,
\[
D_j^{int} = \{ |x| \in [\l(1-2^{-2(j-1)}),\l(1-2^{-2(j+1)})] \}, \qquad 1
\leq 2^j \leq \l^{\frac23}
\]
By $D^{bd}$ we denote a narrow strip near the boundary of the ball,
\[
D^{bd} = \{ ||x|-\l| \leq  \l^{-\frac13}\}
\]
Finally, $D^{ext}$ is the exterior of the ball, 
\[
D^{ext} = \{ |x| > \l + \frac12  \l^{-\frac13} \}
\]

The symbol of $x^2 -\Delta -\l^2$ equals  
\[
x^2+\xi^2 -\l^2
\]
If $x \in D_j^{int}$ then this symbol can be zero only in the 
region 
\[
|\xi| \approx \l 2^{-j}
\]
Hence the case $2^{j} \approx \l^{\frac23}$ corresponds to an annulus
of thickness $\l^{-\frac13}$ and to $|\xi| \approx \l^\frac13$ which
is exactly on the scale of the uncertainty principle. This explains
the limitation for the range of $2^j$. It also gives the scale on
which $\phi$ decays away from the ball.

Given an eigenvalue $\l$ we define the spaces
$l^q_\l L^p$ of functions in $\R^n$ with norm
\[
\| f \|_{ l^q_\l L^p}^q =\| f\|_{L^p(D^{ext})}^q+ \|
f\|_{L^p(D^{bd})}^q +\sum_{1 \leq j, 2^j \leq \lambda^{2/3} } 
\| f\|_{L^p(D_j^{int})}^q
\]
with the usual modification when $q=\infty$. The subscript $\l$ is
used because the sets $D_j^{int}, D^{bd}$ and $D^{ext}$ depend on
$\l$. 

For $x \in \R^n$ we denote
\[
y = \l^{-\frac23}(\l^2-x^2), \qquad \langle y \rangle_- = 1+ y_-,
\qquad \langle y \rangle_{+} = 1+ y_+.
\]
Then

\begin{theorem} 
a) Let $ 2 \leq p \leq  \frac{2(n+1)}{n-1}$. Then
\begin{equation}
\| \l^{\frac13 - \frac{n}3(\frac12-\frac1p)} 
\langle y\rangle_+^{-\frac14 + \frac{n+3}4 (\frac12 -\frac1p)}  
\langle y\rangle_-^{1-\frac{n}2(\frac12-\frac{1}p)}  \phi\|_{l^\infty_\l L^p}
\lesssim   \|\phi\|_{L^2} + \|(H-\l^2)\phi\|_{L^2}
\label{7}\end{equation}

b) Let $ \frac{2(n+1)}{n-1} \leq p \leq \infty$. Then
\begin{equation}
\| \l^{\frac13 - \frac{n}3(\frac12-\frac1p)} 
\langle y\rangle_+^{ \frac12 - \frac{n}2(\frac12 -\frac1p) }
\langle y\rangle_-^{N}  P_{\l^2} \phi\|_{l^\infty_\l L^p}
\lesssim   \|\phi\|_{L^2} 
\label{8}\end{equation}
\label{weight}\end{theorem}

In the interesting regions $D_j$ the estimates can be reformulated as
follows. If $2\le p \le \frac{2(n+1)}{n-1}$ then
\begin{equation}\label{low} 
 \l^{\frac12 - \frac1p} 2^{\frac{j}2( 1- (n+3)(\frac12-\frac1p))}  
\| \phi\|_{L^p(D^{int}_j)}
\lesssim   \|\phi\|_{L^2} + \|(H-\l^2)\phi\|_{L^2}
\end{equation}
and if $ \frac{2(n+1)}{n-1} \leq p \leq \infty$ then
\begin{equation} \label{high} 
 \l^{1 -n(\frac12-\frac1p)} 2^{j(1-n(\frac12-\frac1p))} 
\| P_{\l^2} \phi\|_{L^p(D_j^{int})}
\lesssim   \|\phi\|_{L^2} 
\end{equation}

The estimates in Theorem~\ref{weight} are sharp. We give here a
heuristic motivation, and complement this with a more precise (yet not
as general) construction in the last section of this paper. 

The symbol of $H$ is $x^2+ \xi^2$.   The $\l^2$ eigenfunctions of $H$
are concentrated in the phase space within a neighbourhood of size $1$
of the characteristic set $\{x^2+\xi^2 = \l^2\}$. For a more precise
description, one can consider the Hamilton flow for $H$ restricted to
the same sphere.  It is periodic, and one can can construct
eigenfunctions which are concentrated in the phase space within a
neighbourhood of size $1$ of each such bicharacteristic. To obtain
eigenfunctions which are pointwise larger one may consider
superpositions of such eigenfunctions corresponding to neighbouring
bicharacteristics.  However, such a concentration can only occur on a
smaller set in the physical space; the optimal balance between the
amplitude and the localization region is dictated by the uncertainty
principle.

Consider first the estimate \eqref{7}. Within $D_1^{int}$ the worst
concentration occurs for eigenvalues concentrated
both in space and in frequency in a $1$ neighborhood of a single
bicharacteristic. For an arbitrary $j$ one needs to consider modes which
within $D_j^{int}$ are spatially concentrated within $2^{-\frac{j}{2}}$ of a
bicharacteristic, with a dual concentration in frequency.  In the extreme
case, for $2^j =\l^\frac23$, the concentration occurs as a bump in a
$\l^{-\frac13}$ ball near the circle $|x|=\l$. Except in dimension
$1$, this corresponds to a solution which is uniformly spread in an
ellipsoid of size $\l \times \l^\frac13$.

A different behavior is responsible for estimate \eqref{8}. Here the
concentration corresponds to an equidistribution of energy between all
bicharacteristics through a point $P$.  For $P \in D_j^{int}$ this
leads to concentration in a $2^j \l^{-1}$ ball around the point. If
for instance $P$ is the origin then this could be a spherically
symmetric mode.
\begin{remark}
  We expect that the estimate holds with $l^\infty_\l L^p$ replaced by
  $L^p$ for all $p > 2$.  This is because the concentration cannot
  occur simultaneously in many regions $D_j^{int}$. From
  Corollary~\ref{lpef} we already know that the improvement holds if
  $\frac{2n}{n-2} \leq p \leq \infty$.  On the other hand, such an
  improvement cannot be true for $p=2$.
\end{remark}

If we want to write an estimate without the weights then the sign of
the power of $\langle y \rangle_+$ becomes important:

\begin{cor}
Let $(\phi, \l^2)$ be an eigenfunction, respectively an eigenvalue for the Hermite 
operator. Then
\[
\| \phi\|_{L^p} \lesssim \l^{\rho(p)} \|\phi\|_{L^2}
\]
where for $n \geq 2$ we have (see Figure~\ref{fig})
\[
\rho(p) = \left\{ \begin{array}{cc} 
\ds -(\frac{1}{2} - \frac{1}{p}) &\ds  2 \leq p < \frac{2(n+3)}{n+1} \cr\cr
\ds -\frac13+ \frac{n}3(\frac12-\frac1p) &\ds 
  \frac{2(n+3)}{n+1} < p \leq \frac{2n}{n-2} \cr \cr
\ds -1 + n(\frac{1}2 -\frac1p) &\ds  \frac{2n}{n-2} \leq p \leq \infty
\end{array} \right.
\]
while for $n=1$
\[
\rho(p) = \left\{ \begin{array}{cc}  \ds  -(\frac{1}{2} - \frac{1}{p}) & 2
    \leq p < 4
\cr\cr
\ds -\frac13+ \frac{1}3(\frac12-\frac1p) & 4 < p \leq \infty
\end{array} \right.
\]    

\label{noweight}\end{cor}    

In the first and the third case the worst bound is in $D_1^{int}$, while in
the second case the worst bound is within $D^{bd}$.

{\em Endpoints discussion:} The arguments below yield the $p=\infty$
endpoint even in dimension $n=1,2$. For $n \geq 3$ the case $p =
\frac{2n}{n-2}$ comes from the Strichartz estimates in the previous
section. The $p = \frac{2(n+3)}{n+1}$ endpoint is false for $n=1$,
where the eigenvalues are simple and have an Airy type behavior at the
ends.  However, it is likely be true for $n \geq 2$, because it is
not possible to have concentration at all scales.

\begin{proof}[Proof of Theorem~\ref{weight}: The $L^{2}$ bound.]

This is the  key step in the proof of Theorem~\ref{weight}, because 
the $L^2$ bound is strong enough to provide the localization 
which is needed for the rest of the arguments.
We assume that 
\[
\| \phi \|_{L^2} + \|(H-\l^2)\phi\|_{L^2}= 1
\]
and we shall prove the following bounds:
\begin{eqnarray}
\|\phi\|_{L^2(D_j^{int})} \lesssim 2^{-\frac{j}2}, &\qquad& \|\nabla
\phi\|_{L^2(D_j^{int})} \lesssim \l 2^{-\frac{3j}2}
\nonumber \\
\label{lll2}
\|\phi\|_{L^2(D^{bd})} \lesssim \l^{-\frac{1}3}, &\qquad& \|\nabla
\phi\|_{L^2(D^{bd})} \lesssim 1
\\ \nonumber
\|(x^2 - \l^2) \phi\|_{L^2(D^{ext})} \lesssim \l^\frac13, & \qquad& 
\|(x^2 - \l^2)^{\frac12} \nabla \phi\|_{L^2(D^{ext})} \lesssim \l^\frac13
\end{eqnarray}
We first offer some intuitive justification for these bounds.  The one
in $D^{ext}$ is an elliptic estimate. For the rest, it is convenient
to think of an almost orthogonal basis for the eigenspaces of $H$
which consists of eigenfunctions which are localized in the phase
space within a neighbourhood of size $1$ of null bicharacteristics of
$H-\l^2$. The energy of such eigenvalues is equidistributed along the
corresponding bicharacteristics, which have length $\l$.  Hence in
order to measure what portion of the energy is contained in
$D_j^{int}$ it suffices to measure the length of the bicharacteristic
within $D_j^{int}$, which is at most $2^{-j} \l$. This justifies the
$L^2$ bounds. For the gradient bounds we simply note that on the 
characteristic set of $H-\l^2$ within $D_j$ we have $|\xi| \approx
2^{-j} \l$.

To prove \eqref{lll2} we begin with some simpler estimates. First note that
\[
\langle \phi, (H-\l^2)\phi\rangle = \| \nabla \phi\|_{L^2}^2 + \langle
(x^2 -\l^2) \phi,\phi \rangle
\]
Since $x^2-\l^2 \geq -\l^2$, this easily leads to 
\begin{equation}
\| \nabla \phi\|_{L^2} \lesssim  \l 
\label{e0}\end{equation}
This argument can be easily improved. Choose 
\[
W = ((\l^2 -x^2)^2 + \l^{\frac43})^{-\frac12}
\]
Then an integration by parts yields
\[
\langle W \phi, (H-\l^2)\phi\rangle = \| W^\frac12 \nabla
\phi\|_{L^2}^2 + O(\|\phi\|_{L^2})
\]
which produces another nonsharp estimate,
\begin{equation}
\| ((\l^2 -x^2)^2 + \l^{\frac43})^{-\frac14} \nabla \phi\|_{L^2}
\lesssim 1
\label{e1}\end{equation}
which is slightly better than \eqref{e0} near the sphere $|x|= \l$.
We use this weaker bound to eliminate some error terms in the
estimates which follow.


We will argue in a manner which is similar to the Carleman estimates.
First we introduce a bounded exponential weight which does not change
the estimates, but allows us to replace the operator $H - \l^2$ with
its conjugate with respect to the weight, namely the operator
$H_{a,\l}$ below. The advantage is that the conjugate operator is no
longer selfadjoint; precisely, the gain in the $L^2$ estimates comes
from the positivity of the commutator between its self-adjoint and its
skew-adjoint part. In order to guarantee this the weight needs to be
chosen roughly so that it is convex along the null bicharacteristics
for $H - \l^2$ near the ball $\{|x|=\l\}$.

An alternate approach would be to obtain a Morawetz type estimate
using a suitably chosen multiplier. We do not pursue this as it seems
slightly less precise, and requires more care in the error estimates.

Consider a weight function $a:\R \to \R$ with the following properties:
 
(i) $a$ is nondecreasing and equals $0$ in $(-\infty,-2]$.

(ii) $a$ and its derivatives satisfy the bounds 
\begin{equation}
| \d^k a (y)| \lesssim (1+|y|)^{\frac12-k}, \qquad k \geq 0 
\label{abd}\end{equation}

Let $c$ be a small positive constant which will be chosen later. All
implicit constants in the estimates that follow will be independent of
$c$. Denote
\[
\psi(x) = e^{ - \l^{-\frac23} a(y)} \phi(x), \qquad y = c \l^{-\frac23} (\l^2 - x^2)
\]
Introduce also the conjugated operator
\[
 H_{a,\l}  = e^{- \l^{-\frac23} a(y)} (H-\l^2) e^{ \l^{-\frac23} a(y)} 
\]
Then
\[
 H_{a,\l} \psi = e^{- \l^{-\frac23} a(y)} (H-\l^2) \phi
\]
The exponential weight is bounded since $x^2 \geq 0$ implies $y \leq
\l^{\frac43}$. Hence 
\[
 \|\phi\|_{L^2} + \|(H -\l^2)\phi\|_{L^2} \approx \|\psi\|_{L^2} + \|  H_{a,\l} \psi\|_{L^2} 
\]
On the other hand, given the above properties of $a$ it is easy to verify that 
$\psi$ satisfies \eqref{lll2} if $\phi$ does. Hence we have replaced
$\phi$ and $H-\l^2$ by $\psi$, respectively $H_{a,\l}$.

The conjugated operator $ H_{a,\l}$ is decomposed into a selfadjoint
and a skew-adjoint part,
\[
H_{a,\l} = H_{a,\l}^{re} +  H_{a,\l}^{im}
\]
\[
 H_{a,\l}^{re} = 
 - \Delta + x^2 -\l^2 - 4 c^2 \l^{-\frac83} x^2 |a'(y)|^2 
\]
\[
  H_{a,\l}^{im} = - 2c \l^{-\frac43} (\d x  a'(y) + a'(y) x \d)
\]
where we use a short notation for operators $\d x a^\prime \phi =\sum_{i=1}^n 
\d_{x_i} ( x_i a^\prime \phi)$.    

Ideally we would like the commutator $\{ H_{a,\l}^{re},
H_{a,\l}^{im}\}$ to be positive definite. However, this is too much to
hope for. Instead, it turns out that we can choose the exponential
weight so that the commutator has a positive symbol only near the
characteristic set of $H_{a,\l}$. To compensate for that we introduce
a real correction term $W_\l$ of the form
\[
W_\l = 4  c^2 b(c \l^{-\frac23} (\l^2 - x^2))
\]
where the positive function $b$ is chosen to satisfy the bounds
\begin{equation}
| \d^k b (y)| \lesssim  \langle y \rangle_-  \langle y
\rangle_+^{-\frac32}  \langle y \rangle^{-k}, \qquad k \geq 0
\label{bbd}\end{equation}

Then we compute
\begin{eqnarray*}
\|  H_{a,\l} \psi \|_{L^2}^2 &=& \|  H_{a,\l}^{re} \psi \|_{L^2}^2+
\|  H_{a,\l}^{im} \psi \|_{L^2}^2 + \langle  [H_{a,\l}^{re}, H_{a,\l}^{im}]
\psi,\psi \rangle
\\ &=& \|  (H_{a,\l}^{re} -W_\l) \psi \|_{L^2}^2+
\|  H_{a,\l}^{im} \psi \|_{L^2}^2 + \langle  C
\psi,\psi \rangle
\end{eqnarray*}
where
\[
C =  [H_{a,\l}^{re}, H_{a,\l}^{im}] + H_{a,\l}^{re} W_\l +  W_\l H_{a,\l}^{re}
- W_\l^2
\]
Therefore we obtain 
\begin{equation}
\langle C \psi,\psi\rangle   \le \|  H_{a,\l}  \psi \|_{L^2}^2 
  \lesssim 1
\end{equation}
Our goal is to choose the weights $a$ and $b$ so that $C$ is a
positive operator which controls the norms in \eqref{lll2} for $\psi$.
Hence we need to compute $C$. This can be rather tedious,
and in order to simplify the analysis we make two conventions:

(a) We discard all derivative  terms which can be controlled by
\eqref{e1}.

(b) All the scalar terms which are negligible are incorporated
into a generic term called ``error'', which satisfies
\[
|\text{error}| \lesssim 1+ c^2 \l^\frac23 \langle y \rangle_+^{-1} \langle y \rangle_-^2
\]

As we shall see later on, the error term in (b) is easily controlled
by the main term in \eqref{mainest} below.
For instance the scalar $W_\l^2$ is an error term. We consider
the two remaining terms in $C$. Using the bounds on the derivatives
of $a$ in (ii) above yields
\[
 [H_{a,\l}^{re}, H_{a,\l}^{im}] = 8c( \l^{-\frac43} x^2 a'(y) - 2 c
 \l^{-2}  \d x a''(y) x \d + \l^{-\frac43} \d a'(y) \d)  + \text{error}
\]
On the other hand
\[
W_\l  H^{re}_{a,\l} +  H^{re}_{a,\l} W_\l = -8  c^2 (\d b(y) \d + (x^2 -\l^2) b(y))
+ \text{error}
\]
Summing up the two sets of estimates and using Cauchy-Schwartz we
obtain
\begin{eqnarray*}
1 \gtrsim&& \int \left( 8 c^2 b(y) - 16 c^2
\l^{-2} x^2 a''_-(y)  -  8 c \l^{-\frac43} a'(y)\right)  |\nabla
\phi|^2 \\ &&
+ \left( 8c \l^{-\frac43} x^2 a'(y) - 8 c^2 b(y)(x^2 -\l^2) + \text{error}\right)
|\phi|^2 dx
\end{eqnarray*}
where $a''_-$ is the negative part of $a''$, $a''_- = (|a''|-a'')/2$.

Using the estimates \eqref{e0} and \eqref{e1} we can control the third
derivative term, and also replace the $x^2$ in the second derivative term
by $\l^2$. The $x^2$ in the first scalar term can be replaced by
$\l^2$ modulo a bounded remainder. Thus we obtain
\begin{equation}
1  \gtrsim\int \left( c^2 (b(y) - 2 a''_-(y)) \right)  |\nabla
\phi|^2 + \left( c \l^{\frac23}  (a'(y) - y b(y)) + \text{error}\right)
|\phi|^2 dx
\label{mainest}\end{equation}

To conclude the proof we need to choose the functions $a$ and $b$ 
so that the two coefficients above are nonnegative and sufficiently
large. The key step is summarized in the next Lemma.

\begin{lemma}
  Let $(\e_k)_{k \geq 1}$ be a slowly varying positive sequence with
  the property that $\sum \e_k = 1$. Then there are functions $a$, $b$
  as above so that in addition we have:
\[
b(y) - 2 a''_-(y) \gtrsim \left\{ \begin{array}{cc} 
(1+ |y|) & y < 1 \cr \cr 
\e_k |y|^{-\frac32} & y \in [2^k,2^{k+1}] \end{array}
\right.
\]
\[
a'(y) - yb(y) \gtrsim \left\{ \begin{array}{cc} 
(1+ y^2) & y < 1 \cr \cr 
\e_k |y|^{-\frac12} & y \in [2^k,2^{k+1}] \end{array}
\right.
\]

\label{lab}\end{lemma}

We first show how to conclude the proof of the $L^2$ estimates
\eqref{lll2}.  First we need to control the error term, and this is
where we use the freedom to choose $c$ sufficiently small.
This is because the positive scalar term in \eqref{mainest} has a 
factor of $c$, while the unbounded part of the error has a factor of
$c^2$. Hence in order to control the error it suffices to insure that
\[
\e_k \gtrsim 2^{-\frac{k}2}
\]
which guarantees that
\[
a'(y) - y b(y) \gtrsim |y|^{-1}
\]

The bounds in the lemma suffice to obtain \eqref{lll2} in the regions
$D^{ext}$ and $D^{bd}$. In order to prove \eqref{lll2} in the region
$D_j^{int}$ we need to also choose the $\e_k$'s so that
\[
\e_k \approx 1, \qquad \text { when } 2^{k+2j} \approx \l^\frac43 
\]

\begin{proof}[Proof of Lemma~\ref{lab}]

(i) The range $y < 0$.
Observe first that we can add any positive function to $b$ for $y < 0$
and improve both inequalities. Hence as long as $a'(0) > 0$ and
$a''(0) \geq 0$ we simply need to choose $b$ sufficiently large and
growing like $|y|$ at $-\infty$.

(ii) The range $ y \geq  0$. Here it suffices to choose $a$ so that 
$a' >0$, $a'' (y) < 0$ for large $y$, $a''(0) \geq 0$  and
\[
a'(y) -2 y a''(y) \gtrsim \e_k |y|^{-\frac12} \qquad  y \in [2^k,2^{k+1}] 
\]
Then $b$ can be chosen $b(y) = a''(y) + (2y)^{-1} a'(y)$ for large $y$
and arbitrary between $2a''_-(y) $ and  $ y^{-1} a'(y)$ for small
$y$.

It remains to describe the choice of $a'$. We begin with a simple
choice, namely 
\[
a'_0(y) = (1+y^2)^{-\frac14}
\]
which satisfies all conditions except for a weaker bound
\[
a'_0(y) -2 y a''_0(y) \gtrsim (1+y^2)^{-\frac54}
\]
Then we choose a nonincreasing function $d$ with $ \frac12 \leq d \leq
\frac32$ so that $d'(0) = 0$ and
\[
d'(y) \sim - \e_k y^{-1}   \qquad  y \in [2^k,2^{k+1}] 
\]
Finally, we set 
\[
a'(y) = d(y) a'_0(y)
\]
Then in $ [2^k,2^{k+1}] $ we have
\[
a'(y) -2 y a''(y) = d(y) (a'_0(y) -2 y a''_0(y)) - 2 y a_0'(y) d'(y)
\gtrsim \e_k |y|^{-\frac12}
\]
\end{proof}

\noindent
{ \em Proof of Theorem~\ref{weight}: the $L^{\frac{2(n+1)}{n-1}}$ bound.}
  We split $D_j^{int}$ into balls $D_{j}^{int,k}$ of radius $r_j = \l
  2^{-2j}$ and consider a corresponding partition of unity
\[
1 = \chi^{bd} + \chi^{ext} + \sum_{j,k} \chi_j^k(x)
\]
with
\[
|\nabla \chi_j^k(x)| \lesssim r_j^{-1}, \qquad |\nabla^2 \chi_j^k(x)| \lesssim  
r_j^{-2}
\]
The localized pieces of $\phi$ are 
\[
\phi_{j}^k = \chi_j^k \phi, \qquad \phi^{bd} = \chi^{bd} \phi, \qquad
\phi^{ext} = \chi^{ext} \phi
\]
We consider three cases.
 
{\bf  a)  The interior estimate.}
We claim that the functions $\phi^k_j$  satisfy
\begin{equation}
\sum_k 2^j \|\phi_{j}^k\|_{L^2}^2 + 2^{-j} \|(H-\l^2) \phi_{j}^k\|_{L^2}^2 \lesssim 1 
\label{l2bound}\end{equation}
The first part follows directly from the lemma, for the second we commute
\[
(H-\l^2) \phi_{j}^k = [ -\Delta,  \chi_{j}^k ] \phi = \nabla   \chi_{j}^k \nabla \phi  
+ \nabla^2 \chi_{j}^k \phi
\]
and use the lemma and the bounds on $\chi_j^k$.

To obtain $L^p$ bounds for each of these pieces we use a Strichartz
type estimate which is a special case of  Theorem~2,~\cite{oi}. For
convenience we state it in the following

\begin{lemma}
Let $W$ be a real potential in the unit ball which satisfies $W \sim 1$ and
\[
 |\d^\alpha W| \lesssim 1, \qquad |\alpha| = 1,2
\]
and let $a^{ij}$ be elliptic coefficients of class $C^2$. 
Then for all $\mu > 1$ and $u$ supported in the unit ball we have
\[
\| u\|_{L^{\frac{2(n+1)}{n-1}}} \lesssim \mu^{-\frac{1}{n+1}} ( \mu^\frac12 \|u\|_{L^2} +
\mu^{-\frac12} \|(a^{ij}\d_{i} \d_j  +\mu^2 W) u\|_{L^2})
\]
\label{se1}\end{lemma}
\begin{proof}[Discussion]
 Strictly speaking, in order to apply
Theorem~2,~\cite{oi} we need the additional bounds
\[
 |\d^\alpha W| \lesssim \mu^{\frac{|\alpha|-2}2}, \qquad |\alpha| \geq 2
\]
But such bounds can be easily gained by truncating the potential $W$
in frequency, 
\[
W = W_0 + W_1  \qquad W_0 = \chi(\mu^{-\frac12} D) W, \quad  
W_1 =(1- \chi(\mu^{-\frac12} D))W
\]
Then the low frequency part $W_0$ satisfies the stronger bounds.  The
high frequency part on the other hand satisfies a pointwise bound
$|W_1| \lesssim \mu^{-1}$ therefore it does not affect the size of the
right hand side.
 Another observation is that Theorem~2,~\cite{oi} only gives 
the $L^{\frac{2(n+1)}{n-1}}$ bound for $u$ in the frequency region
$|\xi | \lesssim \mu$. However, outside this region the symbol
$-\xi^2+\mu^2 W$ is elliptic, therefore even stronger bounds are
easy to obtain.
\end{proof}

Given the $L^2$ bound \eqref{l2bound}, the $L^\frac{2(n+1)}{n-1}$
estimate for $\phi$ in $D_j^{int}$ would follow from an estimate for
$\phi_j^k$, namely
\[
\|\phi_j^k\|_{L^{\frac{2(n+1)}{n-1}}} \lesssim (2^j \l^{-1})^{\frac{1}{n+1}}
( 2^j \|\phi_{j}^k\|_{L^2} + 2^{-j} \|(-\Delta + x^2 -\l^2) \phi_{j}^k\|_{L^2})
\]
This follows  from Lemma~\ref{se1} after rescaling to the unit spatial
scale.  Note that within $D_j^{int,k}$ the symbol $x^2+\xi^2 -\l^2$ is
elliptic at frequencies $|\xi| \gg 2^{-j} \l$, so the interesting
region in phase space has size $(2^{-2j} \l)^n \times (2^{-j}\l)^n$.
After rescaling, the frequency becomes $\mu = 2^{-3j} \l^2 > 1$.

{\bf b) The boundary estimate.} 
We can use the bounds in \eqref{lll2} to obtain
\[
\l^\frac13 \| \phi^{bd} \|_{L^2} +  \| \nabla \phi^{bd} \|_{L^2}
\lesssim 1
\]
Then the $L^{\frac{2(n+1)}{n-1}}$ bound for $\phi^{bd}$ is
straightforward by Sobolev embeddings.

{\bf c) The exterior estimate.} 
The bound in \eqref{lll2}
implies that
\[
\| (x^2 -\l^2) \phi^{ext}\|_{L^2} 
+ 
\| (x^2-\l^2)^\frac12 \nabla \phi^{ext}\|_{L^2}
\lesssim \l^\frac13
\]
Then  by (weighted) Sobolev embeddings we obtain 
\[
\| (x^2-\l^2)^\frac{n+2}{2(n+1)} \nabla \phi^{ext}\|_{L^{\frac{2(n+1)}{n-1}}}
\lesssim \l^\frac13
\]
\end{proof}

\begin{proof}[Proof of Theorem~\ref{weight}: the
  $L^\infty$ bound.]

We consider separately two cases.

(i) In $D^{int}$ and $D^{bd}$ we only use the size of the potential.
Each region $D_j^{int,k}$ has size $r_j = 2^{-2j} \l$ and corresponds
to frequencies of size $\mu_j = 2^{-j} \l$.
In $D_j^{int,k}$ we have an elliptic estimate,
\[
\mu_j^{n(\frac1q-\frac1p)} 
\| \chi_j^k \phi\|_{L^q}  \lesssim  \mu_j^{-2} \| \Delta \phi\|_{L^p(D_j^{int,k})}
+ \|\phi\|_{L^p(D_j^{int,k})} , \qquad 1 < p \leq q \leq \frac {2p}{n-2p}
\]
which follows by rescaling from the case $r=1$, $\mu \geq 1$.

We replace $ \Delta \phi$ with $(\Delta -x^2+\l^2) \phi + (x^2-\l^2) \phi$, begin 
with $p = \frac{2(n+1)}{n-1}$ and apply the above estimate iteratively
until we arrive at $q = \infty$. The same idea works also in $D^{bd}$.

(ii) In $D^{ext}$ we also use the favorable sign of the potential.
We first improve the $L^2$ bound on $\phi$, namely the last part 
of \eqref{lll2}. We rewrite \eqref{lll2} as 
\[
\| \langle y\rangle_+ \phi\|_{L^2(D^{ext})} \lesssim \l^{-\frac13}
\]
and inductively show that 
\begin{equation}
\| \langle y\rangle_+^{\frac{N}2} \phi\|_{L^2(D^{ext})} \lesssim \l^{-\frac13}
\label{l2bt}\end{equation}
Suppose that \eqref{l2bt} holds. Then we
compute
\begin{eqnarray*}
\int  \chi^{ext} \langle y\rangle_+^{\frac{N}2} \phi \  \langle
y\rangle_+^{\frac{N}2} (H-\l^2)\phi \, dx \!\! &=& \int  \chi^{ext}
 \langle y\rangle_+^N ((x^2 -\l^2) \phi^2 + |\nabla \phi|^2 )
\\ &+&   \chi^{ext} \langle y\rangle_+^{N-1} \l^{-\frac23}  \phi
x \nabla \phi + \nabla \chi^{ext}   \langle y\rangle_+^N \phi
\nabla \phi dx
\end{eqnarray*}
In the region $y \gg 1$ the last two right hand side terms are
controlled by the first. In the region $y \sim 1$ we use
\eqref{lll2}. Together with \eqref{l2bt} this yields
\[
 \int  \chi^{ext}
 \langle y\rangle_+^N ((x^2 -\l^2) \phi^2 + |\nabla \phi|^2 ) dx
 \lesssim 1
\]
which implies \eqref{l2bt} with $N$ replaced by $N+1$.

At this point we can conclude the proof as in case (i).  Precisely, to
fix the size of the potential we also consider a covering with balls
\[
D^{ext} \subset \bigcup_{j\geq 1,k} D_j^{ext,k}
\]
and a corresponding partition of unit
\[
1 = \sum \chi_j^k \qquad \text{in } D^{ext}
\]
so that  in $D_j^{ext,k}$ we have
\[
|x|^2 -\l^2 \sim 2^{2j} \l^{\frac23}
\]
In each of these balls we use the same elliptic argument as in (i).
\end{proof}

\section{Extensions}

In the previous section we used the potential $x^2$. However, its
precise form does not play a fundamental role in the estimates.  Here
we consider instead the operator $H_V$ with a positive potential $V$
which satisfies the following conditions:
\begin{equation}
 V \sim |x|^2, \qquad |\nabla V| \sim |x|, \qquad |\d_x^2 V|
 \lesssim 1
\label{V} \end{equation}
Given an eigenvalue $\l^2$ of $V$ we introduce as before the dyadic 
regions
\[
D_j^{int} = \{ \l^2 - V \in [ 2^{-2(j-1)}  \l^2,2^{-2(j+1)}  \l^2] \qquad
2 \leq 2^j \leq \l^{\frac23} \}
\]
\[
D^{bd} = \{ |\l^2 - V| \lesssim \l^{\frac23}\} \qquad 
D^{ext} = \{V > \l^2 + \l^{\frac23}\}
\]
Also we set 
\[
y = \l^{-\frac23} (V - \l^2)
\]
Then

\begin{theorem}\label{weighted}  
a) Let $ 2 \leq p \leq  \frac{2(n+1)}{n-1}$. Then
\begin{equation}
\| \l^{\frac13 - \frac{n}3(\frac12-\frac1p)} 
\langle y\rangle_-^{-\frac14 + \frac{n+3}4 (\frac12 -\frac1p)}  
\langle y\rangle_+^{1-\frac{n}2(\frac12-\frac{1}p)}  \phi\|_{l^\infty_\l L^p}
\lesssim   \|\phi\|_{L^2} + \|(H_V-\l^2)\phi\|_{L^2}
\end{equation}

b) Let $ \frac{2(n+1)}{n-1} \leq p \leq \infty$. Then
\begin{equation}
\| \l^{\frac13 - \frac{n}3(\frac12-\frac1p)} 
\langle y\rangle_-^{ \frac12 - \frac{n}2(\frac12 -\frac1p) }
\langle y\rangle_+^{N}  P_{[\l^2,\l^2+1]} \phi\|_{l^\infty_\l L^p}
\lesssim   \|\phi\|_{L^2} 
\end{equation} 
\end{theorem}

The result in Corollary~\ref{noweight} also applies to $H_V$.

\begin{proof}(Sketch)
We only need the counterpart of \eqref{lll2}, the rest of the
proof is identical. 

Without any restriction in generality we can
replace $V$ with a mollified potential $\chi(D) V$, since the
difference is bounded. This allows us to use the additional assumptions
\eqref{vi}.

We use the same functions $a$ and $b$ as in the proof of \eqref{lll2} but
with the above definition of $y$.  Now the selfadjoint and
skew-adjoint parts of the conjugated operator are
\[
H^{re}_{a,\l} = -\Delta + V- \l^2 - c^2 \l^{-\frac63} |\nabla V|^2 |a'(y)|^2
\]
\[
H^{im}_{a,\l} = - c \l^{-\frac43} (\d \nabla V a'(y) + a'(y) \nabla V \d)
\]
The correction term is chosen to be 
\[
W_\l = 2 c^2 \l^{-2} |\nabla V|^2 b(y)
\]
Modulo error terms the operator $C$ has the form
\begin{eqnarray*}
C &=& 2c \left(\l^{-\frac43} |\nabla V|^2 a'(y)  -  \l^{-2}  b(y) |\nabla
V|^2 (V-\l^2)\right) + 2 c^2 \d( \l^{-2} |\nabla V|^2 b(y) \\ &+&
\l^{-2} \nabla V a''(y) \nabla V) \d
\end{eqnarray*}
Then the argument continues  as in the proof of \eqref{lll2}.

\end{proof}

\section{Remarks on optimality}

We begin by recalling bounds and expansions of Hermite functions.
The Hermite functions are the eigenfunctions of the one dimensional
Hermite operator, and solve
\begin{equation}\label{hermite}
 -h^{\prime \prime }_k+ x^2 h_k = (2k+1) h_k  
\end{equation} 
for nonnegative integers $k$. 
They are given by
\[  
h_k(x) = e^{x^2/2} (-1)^k \frac{d^k}{dx^k} e^{-x^2}.  
\]
and are even functions for even $k$ and odd functions for odd $k$.
Here the meaning of $k$ differs slightly from the previous sections.  
In dimension $n$ a complete set of eigenfunctions is given by
\[
h_\alpha (x) = \prod_{i=1}^n  h_{\alpha_i}(x_i)
\]
where the corresponding eigenvalue is $n+ 2|\alpha|$.

To construct highly localized eigenfunctions we need a better
understanding of the behavior of the Hermite functions. This is well
understood by now, and we describe it next.

The ODE \eqref{hermite} has a turning point at $x= \sqrt{2k+1}$.  We
set $\lambda = \sqrt{2k+1}$. Then the functions $h_k$ have an
oscillatory behavior for small $x$, an Airy type behavior for $|x|$
close to $\l$ and Gaussian decay  for large $x$.  More precisely,
define
\begin{equation} 
 s^-(x) = \int_{0}^x     \sqrt{|t^2 -\lambda^2|} dt,  
\end{equation} 
\begin{equation} 
 s^+(x) = \int_{\lambda}^x     \sqrt{|t^2 -\lambda^2|} dt,  
\end{equation} 
Then 
\begin{lemma}
 \label{hermiteest}
The normalized eigenfunctions 
\[
\tilde{h}_k = h_k \|h_k\|_{L^2}^{-1}
\]
satisfy
\[
\tilde{h}_{2k} = \left\{ \begin{array}{lc} a_{2k}^- 
 (\l^2-x^2)^{-\frac14}(\cos{s^-(x)} +
  error) & |x|<\lambda -\l^{-\frac13} \cr O(\l^\frac16) & \lambda
  -\l^{-\frac13} \leq x \leq  \lambda+\l^{-\frac13} \cr
 a_{2k}^+ e^{- s^+(x)}  (\l^2-x^2)^{-\frac14}(1 +
  error) & |x|> \lambda+\l^{-\frac13} \end{array} \right.
\]
\[
\tilde{h}_{2k+1} = \left\{ \begin{array}{lc} a_{2k+1}^-  (\l^2-x^2)^{-\frac14}(\sin{s^-(x)} +
  error) & |x|<\lambda - \l^{-\frac13} \cr O(\l^\frac16) & \lambda
  -\l^{-\frac13} \leq x \leq  \lambda+\l^{-\frac13} \cr
 a_{2k+1}^+ e^{- s^+(x)}  (\l^2-x^2)^{-\frac14}(1 +
  error) & x >\lambda + \l^{-\frac13}\end{array} \right.
\]
where
\[
|a_k^\pm| \sim 1, \qquad error = O((|x^2-\l^2|^{-\frac12}
||x|-\l|^{-1})
\]
\end{lemma}
Note that the error term is $O(1)$ if $||x|-\l| \sim \l^{-\frac13}$
and decays away from $\l$. One can also write an Airy type asymptotic 
near $|x| = \l$, but we do not need it here.
The bounds of Lemma \ref{hermiteest} follows from standard WKB bounds as in 
\cite{MR95m:34091}
and well-known formulas for the Hermite functions, see \cite{MR2000g:33001}.


An immediate consequence of the above lemma is that all bounds in
Theorem \ref{weighted} are sharp in one dimension.  In what follows we
construct examples which show that the bounds in Theorem
\ref{weighted} are also sharp in higher dimension. For a positive
integer $N$  we consider eigenfunctions which correspond to the 
eigenvalue $n+2N$. We define $\l > 0$ by
\[ 
\lambda^2 = n + 2N 
\]

\subsection{ Concentration in a tube in $D_0$.}
 We set $x^\prime = (x_2, \dots, x_n)$ and consider the eigenfunction
\[ 
v(x)  = \tilde h_N(x_1) e^{-\frac12|x^\prime|^2 } 
\] 
which is concentrated in the tube 
\[
T = \{ |x_1| < \l, \qquad |x'| < 1\}
\]
We have
\[ 
\Vert v \Vert_{L^2} \sim 1 
\]
and, if $1\le p \le \infty$  
\[ \Vert v \Vert_{L^p(T \cap D_0^{int})} \sim \lambda^{\frac1p-\frac12}. \]

\subsection{ Point concentration  in $D_0$.}
 We assume that $N$ is even, otherwise the arguments need a small 
modification. We consider the set $I$ of indices
\[
I = \{ \alpha\  \text{even} ; |\alpha|=N,\ \alpha_i > N/4n \}
\]
Then we consider the eigenfunction
\[
v = \sum_{\alpha \in I } \prod_{i=1}^n h_{\alpha_i}(x_i)
\]
which concentrates in the ball $B(0,\l^{-1})$.
Since 
\[
|I| \sim  N^{n-1} \sim  \lambda^{2( n-1)}
\]
it follows that  
\[
\|v\|_{L^2} \sim  \lambda^{n-1}
\]
On the other hand, by Lemma~\ref{hermiteest}, for $\alpha \in I$ we get  
\[
h_{\alpha_i}(x_i) \sim \l^{-\frac12}, \qquad |x_i| < \l^{-1} 
\]
Summing up, we obtain
\[
h(x) \sim |I| \l^{-\frac{n}2} \sim  \lambda^{2( n-1)}
\lambda^{2( n-1)}, \qquad x \in B(0,\l^{-1})
\]
therefore
\[
 \Vert h \Vert_{L^p(B(0,\l^{-1}))} \sim \lambda^{-\frac{n}p}\lambda^{-n/2} 
 \lambda^{2(n-1)} 
\]
Thus 
\[ 
\lambda^{1-n(\frac12-\frac1p)} \Vert v   \Vert_{L^p(D^{int}_0)} \gtrsim
\Vert v \Vert_{L^2} 
\]
which one should compare to \eqref{high}.

\subsection{Point concentration in $D_j$}
We fix $r$ with $re_1 \in D^{int}_j$ and construct an eigenfunction
which concentrates in $B(r e_1, c 2^j\lambda^{-1})$. Let $N_0$ 
be so that
\[
|N_0 - \frac{\l^2 -r^2}2 | < 1.
\]
This implies that
\[
|N_0| \sim 2^{-2j} \l^2
\]

We consider the set of indices 
\[
I = \{ \alpha| \ \alpha'\  \text{even} ; |\alpha|=N,\ \alpha_1 > N
-N_0,\ \alpha'_i > N/4n \}
\]
and its subset
\[
J = \{ \alpha \in I; \ h_{\alpha_1}(r) > \frac14 (\l^2-r^2)^{-\frac12}\}
\]
Lemma~\ref{hermiteest} shows that near $r$ the functions
$h_{\alpha_1}$ oscillate at a frequency of the order of $2^{-j} \l$.
Then, changing $r$ by no more that $2^j\lambda^{-1}$, we 
can insure that $|J| \sim |I|$. We define
\[
v =  \sum_{\alpha \in J } \prod_{i=1}^n h_{\alpha_i}(x_i) 
\]
 
Since 
\[
|I| \sim N_0^{n-1}
\]
we have 
\[
\|v\|_{L^2} \sim |J|^\frac12 \sim 2^{-(n-1)j} \l^{n-1}
\]
On the other hand, using Lemma~\ref{hermiteest} and the definition 
of $J$ we compute
\[
v(x) \sim 2^{\frac{j}2} \l^{-\frac12} N_0^{-\frac{n-1}4} |J| 
\sim \l^{-\frac{n}2} 2^{\frac{n-1}2 j} 2^{-2(n-1)j} \l^{2(n-1)},
\qquad x \in B(r e_1, c 2^j\lambda^{-1})
\]
This gives
\[
\| v\|_{L^p(B(r e_1, c 2^j\lambda^{-1}))} \sim 2^{\frac{nj}p} \lambda^{-\frac{n}p} 
 \l^{-\frac{n}2} 2^{\frac{n}2 j} 2^{-2(n-1)j} \l^{2(n-1)}
\]
therefore
\[
\lambda^{1-n(\frac12-\frac1p)} 2^ {j(1- n(\frac12-\frac{1}p) } \Vert v
\Vert_{L^p(D^{int}_j)} \gtrsim \vert v \Vert_{L^2}.
\]

\subsection{Concentration on a tube in $D_j^{int}$.} Here we construct
eigenfunctions which concentrate in the tube
\[
T = \{ x \in D_j^{int};\ |x'| \leq 2^{-\frac{j}2} \}
\]
Now we consider the set of indices 
\[
I = \{ \alpha; \ |\alpha|=N,\ \alpha' \text{ even},\  |\alpha'_j -
2^j| \leq c 2^j \}
\]
with small $C$. This has size
\[
|I| \sim 2^{j(n-1)}
\]
We want to consider a subset $J$ of $I$ of comparable size 
and set 
\[
v =  \sum_{\alpha \in J } \prod_{i=1}^n h_{\alpha_i}(x_i) 
\]
The difficulty is that we want to avoid cancellations in this
summation.  For $\alpha \in I$ the functions $h_{\alpha'_i}(x'_i)$ are
of size $2^{-j/4} $ in an interval of length $2^{-j/2}$ centered at
zero. It remains to insure that the functions $h_{\alpha_1}(x_1)$
have the same sign within most of the tube. For this we need to check 
that their phases are coherent, i.e. they differ essentially by a constant.

Set $\mu = \sqrt{1+\alpha_1}$. We must compare the functions
$s_\mu^-(x)$ in the range
\[
\l -|x| \sim 2^{-2j}\l, \qquad  \l^2 - \mu^2 \sim 2^j 
\]
For this we compute
\[
\frac{d}{d\mu} \frac{d}{dx} s_\mu^-(x) = \frac{2\l}{\sqrt{\l^2-x^2}}
\sim 2^j
\]
Integrating this with respect to $x$ and then with respect to $\mu$ we find 
that for $x$ and $\mu_1$, $\mu_2$ as above we have
\[
s_{\mu_1}^-(x) - s_{\mu_2}^-(x) = C(\mu_1,\mu_2) + O(2^{-j} \l |\mu_1-\mu_2|)
\]
which yields
\[
|s_{\mu_1}^-(x) - s_{\mu_2}^-(x) - C(\mu_1,\mu_2)| \lesssim c
\]
If $c$ is small enough then the phases $s_{\mu}^-(x)$ are close modulo
constants, therefore we can choose a subset of indices $\mu$ of
comparable cardinality so that the phases are close modulo $2 \pi \Z$.
This leads to the subset $J$ of $I$ corresponding to this restricted
set of $\alpha_1$'s.

Now we can complete the computation.
On one hand
\[
\|v\|_{L^2} \sim 2^{\frac{n-1}2 j}
\]
On the other hand, for most $x \in T$ we have
\[
|v(x)| \sim 2^{\frac{j}2} \l^{-\frac12} 2^{-\frac{n-1}{4}j} 
|J| \sim  2^{(n-1)j} \l^{-\frac12} 2^{-\frac{n-3}{4}j} 
\]
Since $|T| \sim 2^{-\frac{n+3}{2}j} \l$ this yields
\[
\|v\|_{L^p(T)} \sim   2^{-\frac{n+3}{2p}j} \l^{\frac{1}p}
2^{(n-1)j} \l^{-\frac12} 2^{-\frac{n-3}{4}j} 
\]
Hence
\[ 
\Vert u \Vert_{L^2} \lesssim  \lambda^{\frac12-\frac1p}
 2^{-\frac{j}2  \frac{n+1}2 } 2^{ \frac{j}2 \frac{n+3}p} \Vert u 
\Vert_{L^p(D^{int}_j)}. 
\]

\end{document}